\documentclass[11pt]{article}
\input epsf
\usepackage{amssymb,amsmath}          
\usepackage[T1]{fontenc}
\usepackage[latin2]{inputenc}
\usepackage{graphicx}

\newtheorem{definition}{Definition}[section]
\newtheorem{prop}{Proposition}[section]
\newtheorem{tm}{Theorem}[section]
\newtheorem{lm}{Lemma}[section]
\newtheorem{pr}{Example}[section]

\title{Testing equality in distribution of random convex compact sets via theory of~$\mathfrak{N}$-distances and its application to assessing similarity of general random sets}
\author{
Vesna Gotovac $^1$, Kate\v{r}ina Helisov\'a$^2$\\
 \\
$^1$Department of Mathematics\\ Faculty of Science\\ University of Split, Croatia\\ 
email: vgotovac@pmfst.hr\\
 \\
$^2$Department of Mathematics\\ Faculty of Electrical Engineering\\ Czech Technical University in Prague, Czech Republic\\
email: helisova@math.feld.cvut.cz\\
 }

\date{December 12, 2017}

\begin{document}

\maketitle

\begin{abstract}
This paper concerns a method of testing equality of distribution of random convex compact sets and the way how to use the test to distinguish between two realisations of general random sets.
The family of metrics on the space of distributions of random convex compact sets is constructed using the theory of $\mathfrak{N}$-distances and characteristic functions of random convex compact sets. 
Further, the approximation of the metrics through its finite dimensional counterparts is proposed, which lead to a new statistical test for testing equality in distribution of two random convex compact sets.
Then, it is described how to approximate a realisation of a general random set by a union of convex compact sets, and it is shown how to determine whether two realisations of general random sets come from the same process using the constructed test.
The procedure is justified by an extensive simulation study.
\end{abstract}


\section{Introduction}
\label{sec:intro}

In the last years, planar random sets have been studied from different points of view because of their widespread applications in biology, material sciences, medicine etc. 
They can describe and explain many events, e.g. behaviour of cells in organisms (see \cite{MM11}, \cite{Hetal15}), particles in materials (see \cite{H14}, \cite{Netal16}) or presence of different plants (see \cite{D81}, \cite{MH10}).
Therefore it is useful to develop methods for their statistical analyses. 

Although it is often beneficial to know the concrete model, 
there are situations when the knowledge about the underlying process is not necessary, for example when we want only to compare random sets based on their realisations like to distinguish between two types of cells based on microscopic pictures, recognise different tendency of the growth of some plants or trees etc. 
Basically, it is similar to classical hypothesis testing, where we can use rank based statistics without specifying strict probability distributions.

The main aim of this paper is to find a statistical procedure for comparing probability distributions of two random sets based of comparing their realisations without the knowledge of the distributions themselves.
A procedure is suggested in \cite{Getal16}. 
It is based on dividing realisations of random sets to unions of convex compact sets using Voronoi tessellation, whose cells are convex compact (see \cite{Cetal13}), and consequent comparing the support functions of the convex sets using envelope test (see \cite{Metal16}), which deals with average values of the support functions.
In the present paper, we propose a new statistical test,
which does not work with the averages only, but with the distributions themselves.
We again focus on the convex compact sets,
because, as mentioned above, we can convert the problem of comparing general sets 
to the problem of testing equality of distributions of convex compact sets, which
have more convenient properties than rugged shaped sets.
Moreover, the analysis of the convex compact sets themselves is very useful since they play and important roles in many models.
They appear for example as grains in popular germ-grain model (see \cite{Cetal13}).
Also the Voronoi diagrams mentioned above have many applications, for example 
as a model of different biological structures like cells or bone micro-architecture. 
In medical diagnosis, Voronoi diagrams-based models for muscle tissue are used to detect neuromuscular diseases (see \cite{Setal16}) etc.

In the present paper, the equality in distribution of two random convex compact sets is tested using non-parametric, distribution-free test based on $\mathfrak{N}$-distances (also called kernel tests).
Two sample tests in finite dimensional spaces constructed by using kernels could be found in \cite{K06} or \cite{Getal12}. 
Our research is focused on construction of 
the kernel on infinite dimensional space whose finite dimensional counterparts could be easily evaluated using available information.
Here, the characteristic function of random convex compact set plays an important role. 
The infinite dimensional analogue of the concept of the characteristic function, named characteristic functional, was first introduced in \cite{K35} for the case of distributions in Banach
space. 
In some particular cases, the domain of this characteristic functional could be replaced with more simple one.
Characteristic functions of random convex compact sets with more simple domain is defined in \cite{L00}, however for our purposes, we slightly modified the definition. 
The important property is that this characteristic function connects the distribution of the random convex compact set with the finite-dimensional distributions of its support function. 
Just note that hypothesis testing using difference of empirical characteristics functions and $L^2$ norm appears also in other fields, e.g. in time series analysis (see \cite{H99}). 

The present paper is organized as follows. 
Section 2 summarizes already existing theoretical results concerning 
random convex compact sets and their characteristic functions, and basics of $\mathfrak{N}$-distance theory. 
In Section 3, there is introduced the characteristic function of random convex compact sets.  
Section 4 is dedicated to testing the equality of their distributions.
In Section 5, the procedure how to use the test for assessing similarity of general random sets is described, and results of the simulation study are shown.

\section{Theoretical background}

In this section, we provide an overview of existing theory and results that we later use to build an infinite-dimensional kernel and related tests. 

\subsection{Convex compact sets}
The definitions and results from this subsection can be found in \cite{Cetal13} and \cite{S13}.
\begin{definition} 
The set $ A \subset \mathbb R^d$ is said to be \textit{convex} if
$cx+(1-c)y\in  A$ for all $x,y \in  A$ and all $0<c<1$.
The set $ A \subset \mathbb R^d$ is said to be \textit{compact} if
it is closed and bounded.
\end{definition}
 By $\mathcal{K}$ we denote family of all convex compact sets in $\mathbb{R}^d.$
\label{sec:conv}
\begin{definition} 
\textit{Support function} of a convex compact set $A$ is a function $h(A,\cdot):\mathbb{R}^d \to \mathbb{R}$  defined by
\[
h(A,u) = \sup_{x \in A} \langle u,x \rangle, \ u  \in \mathbb{R}^d.
\]
\end{definition}
Since $\langle u, \cdot \rangle$ is a continuous function and the set $A$ is compact, $h(A,u)$ is well defined and supremum is attained for all $u \in \mathbb{R}^d$. 
Also, it is easy to verify that support function is convex and therefore it is continuous.
Moreover since $h(A,\lambda u)=\lambda h(A, u),$ the support function $h(A, \cdot)$ is uniquely determined by its values on $S^{d-1},$ where $S^{d-1}$ is the unit sphere in $\mathbb{R}^d.$ Therefore, we focus on the case when the domain of $h(A, \cdot)$ is $S^{d-1},$ while we consider the variables $u$ to be the corresponding angles instead of points in $\mathbb{R}^d.$

\begin{figure}[t]
\centering
\includegraphics[width=15cm]{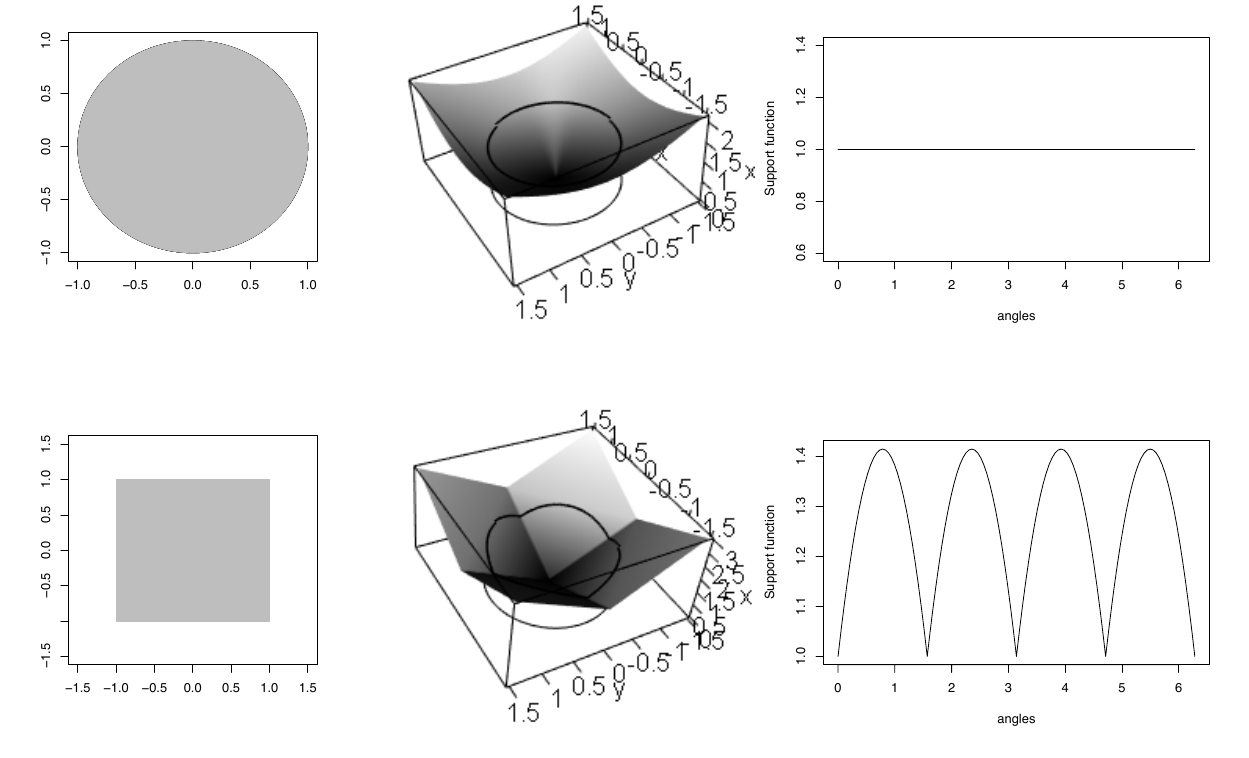}
\caption[Examples of the support function for a disc and a square] {Examples of the support function for a disc and a square,
respectively, when the domain is $\mathbb{R}^2$ (middle) and when the domain is $S^1$ (right).}
\label{fig:supportexample}
\end{figure}

For $u \in S^{d-1},$ the value of $h(A,u)$ represents the signed distance from the origin to the support plane $H(A,u).$ This distance is positive if both $A$ and the origin lie on the same side of the support plane, negative if $A$ and the origin are separated by the support plane and equal to zero if the origin lies in the support plane.
Examples of the support functions of a disc and a square are given in Figure \ref{fig:supportexample}.

\begin{prop}
\label{kor:intersec}
Each non-empty closed convex set in $\mathbb{R}^d$ is the intersection of all its supporting half-spaces $H^{-}(A,u):=\left\{x \in\mathbb{R}^d | \langle x,u\rangle \leq h(A,u) \right\}, u \in S^{d-1}.$
\end{prop}

As follows from Proposition $\ref{kor:intersec},$ each convex compact set is uniquely determined by its support function.

According to the map $\varphi : \mathcal{K} \rightarrow C(S^{d-1}),$
$\varphi(A) = h(A,\cdot),$
$\mathcal{K}$ is embedded into a space $C(S^{d-1})$ of continuous real functions on $S^{d-1}.$

\subsection{Random convex compact set}
\label{sec:r_convex}
Let $\mathcal{F}$ be the family of closed sets and $\mathcal{C}$ the family of compact set of topological space $\mathbb{R}^d$ with the standard topology $\mathcal{G}$.

The \textit{Effros} $\sigma$-algebra $\mathcal{B}(\mathcal{F})$ on $\mathcal{F}$ is the $\sigma$-algebra generated by sets $\left\{ F \in \mathcal{F}: F \cap C \neq 0 \right\}$ for all $C \in \mathcal{C}.$

\begin{definition} 
\label{def:randomset}
Let $(\Omega, \Sigma, P)$ be a probability space.
Then a \textit{random closed set} $\mathbf X$ in $\mathbb R^d$ is a measurable
mapping from $(\Omega, \Sigma)$ to $(\mathcal F,\mathcal B (\mathcal F))$.
A random closed set $\mathbf X$ in $\mathbb{R}^d$ with 
almost surely convex compact values is called \textit{random convex compact set.}
\textit{The distribution $P_{\mathbf X}$ of a random set} \ $\mathbf X$ \ is \ given \ by \ the \ relation $P_{\mathbf X}(F) =
P\left( \left\{
  \omega \in \Omega: \mathbf X(\omega) \in F
\right\} \right)$
for $F \in \mathcal{B}(\mathcal{F})$.
\end{definition}
Denoting $d_E$ the Euclidean metric on $\mathbb{R}^d,$ we can consider Hausdorff metric on $\mathcal{F}$ defined by
\begin{equation}
d_H(A,B)=\max\left\{ \sup_{x \in A}\inf_{y \in B}d_E(x,y), \sup_{y \in B}\inf_{x \in A}d_E(x,y)\right\}, \ A,B \in \mathcal{F}.
\end{equation}
For $A,B \in \mathcal{K}$ it holds that
\begin{equation}
d_H(A,B)=\sup\limits_{u \in S^{d-1}}\left|h(A,u)-h(B,u) \right|,
\label{eq:haussdorf_char_support}
\end{equation}
see \cite{S13}.
\begin{tm}
\label{tm:Effros_Hausdorff}
The Effros $\sigma$-algebra induced on $\mathcal{K}$ coincides with the Borel $\sigma$-algebra generated by Hausdorff metric.
\end{tm}
For the proof, see \cite{S13}.

Further on, let us denote by $\mathcal{U}$ Borel $\sigma$-algebra on $\mathcal{K}$ generated by Hausdorff metric. 
\begin{tm}
\label{tm:char_distr}
For two random convex compact sets $\mathbf X_1$ and $\mathbf X_2$, it holds that
\[
  \mathbf X_1 =^{(\mathcal{D})} \mathbf X_2 
  \Leftrightarrow 
  \left(
    h(\mathbf X_1,u_i)
  \right)_{i\in I}
  =^{(\mathcal{D})}
  \left(
    h(\mathbf X_2,u_i)
  \right)_{i\in I}
\]
for all finite index sets $I$ and $(u_n)_{n \in \mathbb{N}}$ dense subset of $S^{d-1}.$
\end{tm}
\noindent
For the proof of this theorem, reader is referred to \cite{L00}.

\subsection{Random continuous functions on $S^{d-1}$ and associated $\sigma$-algebras}

Since the space of the support functions of the convex compact sets forms a subset of $C(S^{d-1}),$ it is useful to focus on random continuous functions on $S^{d-1}$ as a space where all the support function of the random convex compact sets lie. 
In this context, it is natural to considered $C(S^{d-1})$  as a metric space with supremum metrics
$d_H(h_1,h_2)=\sup\limits_{u \in S^{d-1}}|h_1(u)-h_2(u)|.$ 
Denote $\mathcal{B}(C(S^{d-1}))$ the Borel  $\sigma$-algebra generated by topology induced by $d_h.$
\begin{definition} 

Let $(\Omega,\Sigma,P)$ be a probability space. A \textit{random continuous function $\mathbf X$ on $S^{d-1}$}, $\mathbf X:\Omega\to C(S^{d-1})$, is a measurable mapping from $(\Omega, \Sigma)$ to $(C(S^{d-1}),\mathcal{B}(C(S^{d-1}))$.
\end{definition}
Note that from Theorem $\ref{tm:Effros_Hausdorff}$ and  relation $(\ref{eq:haussdorf_char_support})$, it follows that according to the map $\varphi : \mathcal{K} \rightarrow C(S^{d-1}),$ $\varphi(A) = h(A,\cdot),$
$(\mathcal{K}, \mathcal{B}(\mathcal{F})\cap\mathcal{K})$ is isomorphically embedded into a space $(C(S^{d-1}),\mathcal{B}(C(S^{d-1})))$, 
so we can identify the notion of random convex compact sets with the notion of its random support function viewed as a random continuous function on $S^{d-1}.$   

Let us further investigate the $\sigma$-algebra $\mathcal{B}(C(S^{d-1})).$ 
We define $\mathcal{B}$ as the smallest $\sigma$-algebra containing all measurable cylinder sets which are the sets of the form
\begin{equation}
\left\{A \in \mathcal{K}: \left(h(u_{1}),\ldots,h(u_{k})\right) \in B_1 \times \ldots \times B_k \right\}
\label{eq:dense_gen}
\end{equation}
for arbitrary $u_1,\ldots,u_k \in S^{d-1}$ and $B_1,\ldots,B_k \in \mathcal{B}(\mathbb{R}).$
Thus, this $\sigma$-algebra can be also called the $\sigma$-algebra of finite-dimensional distributions as the cylinder sets exactly correspond to the distributions of support function evaluated at the finitely many angles.

\begin{prop}
For any $D$ dense countable subset of $S^{d-1}$, $\sigma$-algebra $\mathcal{B}$ is generated 
by the family of the sets of the form \eqref{eq:dense_gen}
for arbitrary $u_1,\ldots,u_k \in D$ and $B_1,\ldots,B_k \in \mathcal{B}(\mathbb{R}).$
\label{prop:dense}
\end{prop}

\begin{prop}
It holds that $\mathcal{B}(C(S^{d-1}))=\mathcal{B}.$
\end{prop}

For the proofs see e.g. \cite{GS96}.

\subsection{Fourier-Stieltjes transform of finite signed measure on $(\mathbb{R}^d,\mathcal{B}(\mathbb{R}^d))$}
\label{sec:four}

In this section, we present the notion of the Fourier-Stieltjes transform of a finite signed measure on $(\mathbb{R}^{d},\mathcal{B}(\mathbb{R}^d)).$  
\begin{definition} 
Let $\nu$ be a finite signed measure on $(\mathbb{R}^d,\mathcal{B}(\mathbb{R}^d)).$ The \textit{Fourier-Stieltjes transform of $\nu$} is the function
$\hat{\nu} : \mathbb{R}^d \to \mathbb{C}$ defined by
\[
\hat{\nu}(t)=\int\limits_{\mathbb{R}^d}\exp{\left(i\left\langle t,x\right\rangle\right)}\nu(dx).
\]
\end{definition}
\begin{tm}[Uniqueness theorem]
Let $\nu_1$ and $\nu_2$ be finite signed measures on $(\mathbb{R}^d,\mathcal{B}(\mathbb{R}^d)).$ It holds that
\[
\nu_1=\nu_2 \Longleftrightarrow \hat{\nu_1}=\hat{\nu_2}.
\]
\label{tm:unique}
\end{tm}
For the proof, see \cite{C75}.

\subsection{Basics of $\mathfrak{N}$-distances}
\label{sec:Ndist_b}
In this subsection, we present a class of distances on the space of the distributions of random elements on arbitrary non-empty set $\mathcal{X}.$ The class is constructed using positive or negative definite kernels. We call them $\mathfrak{N}$-distances. For all definitions, propositions and proofs in this subsection, see \cite{K06} or \cite{Retal13}.

\begin{definition} 
Let $\mathcal X$ be a non-empty set. A map $ K: \mathcal X \times \mathcal X \rightarrow \mathbf C$ is 
called \textit{positive definite kernel} if for any $n \in \mathbf N$, arbitrary $c_1,..., c_n 
\in \mathbf C$ and arbitrary $x_1,..., x_n \in \mathcal X$ it holds that
$$
\sum_{i=1}^n \sum_{j=1}^n K(x_i, x_j)c_i \bar c_j \geq 0.
$$
\end{definition}
\begin{definition} 
Let $\mathcal X$ be a nonempty set. A map $ L: \mathcal X \times \mathcal X \rightarrow \mathbf C$ is 
called \textit{negative definite kernel} if for any $n \in \mathbf N$, arbitrary $c_1,..., c_n \in \mathbf C$ 
such that $\sum_{i=1}^n c_i = 0$ and arbitrary $x_1,..., x_n \in \mathcal X$ it holds that
\begin{equation}
\sum_{i=1}^n \sum_{j=1}^n  L(x_i, x_j)c_i \bar c_j \leq 0.
\label{eq:ndk}
\end{equation}
\end{definition}

\begin{prop} 
A function $\mathcal K: \mathcal X \times \mathcal X \rightarrow \mathbf C$ is positive definite kernel if 
and only if there exists a family $\{ f_i \}_{i\in I}$ of complex-valued functions such that
 $\sum_i |f_i(x)|^2 < \infty$ for all $x \in \mathcal X$ and
$$
 K(x,y) = \sum_i f_i(x) \bar f_i(y), \ x,y \in \mathcal X.
$$
\label{prop:pdk}
\end{prop}
\begin{prop}
If $ K: \mathcal X \times \mathcal X \rightarrow \mathbf C$ is a positive definite kernel, then
$$
 L(x,y) =  K(x,x) +  K(y,y) - 2 K(x,y), \ x,y \in \mathcal X
$$
is negative definite kernel such that $ L(x,y) = \overline{ L(y,x)}$ and $L(x,x)=0$.
\label{prop:ndk}
\end{prop}

In the sequel, we suppose that $\mathcal{X}$ is a metric space. 
We denote $\mathfrak{U}$ the Borel $\sigma$-algebra with respect to the topology induced by the metric. 
Consider negative definite kernels on $\mathcal{X}$ which are continuous, symmetric and real-valued. Denote by $\mathcal{P}$ the set of all probability measures on $(\mathcal{X},\mathfrak{U}).$
Suppose that $L$ is a real continuous function, and denote $\mathcal{P}_{L}$ the set of all measures $\mu \in \mathcal{P}$ for which the integral
$
\int_{\mathcal{X}}\int_{\mathcal{X}}L(x,y)d\mu(x)d\mu(y)
$
exists.
\begin{prop}
Let $L$ be a real continuous function on $\mathcal{X} \times \mathcal{X}$ under the condition 
$L(x,y)=L(y,x), \ x,y \in \mathcal{X}.$
The inequality
\begin{align}
2\int_{\mathcal X} \int_{\mathcal X} &  L(x,y) d\mu (x) d\nu (y) - \int_{\mathcal X} \int_{\mathcal X}  L(x,y) d\mu (x) d\mu (y)
-\int_{\mathcal X} \int_{\mathcal X}  L(x,y) d\nu (x) d\nu (y)\geq 0 
\label{eq:ndist}
\end{align}
 holds for all $\mu, \nu \in \mathcal{P}_{L} $ if and only if $ L$ is a  negative definite kernel.
\end{prop}
\begin{definition} 
Let $Q$ be a probability measure on $\mathcal X$ and $c$ be a function on $\mathcal X$ such that
$\int_{\mathcal X} c(x) dQ(x)=0.$ Then $ L$ is called \textit{strongly negative definite kernel} if for any measure $Q,$
$$
\int_{\mathcal X} \int_{\mathcal X}  L(x,y) c(x) c(y) dQ(x) dQ(y) = 0
$$
implies $c(x)=0$ $Q$-almost everywhere.
\label{df:strongly}
\end{definition}
\begin{prop}
Let $L(x,y) =  L(y,x)$. Then \eqref{eq:ndist}
holds for all measures $\mu, \nu \in \mathcal P_{\mathcal L}$ with 
equality in the case $\mu=\nu$ only, if and only if $ L$ is a strongly negative definite kernel.
\label{prop:strong}
\end{prop}
Let us now give few examples of the strongly negative kernels.

\begin{pr}
\label{ex:R_d_kernels}
Strongly negative definite kernels on $\mathbb{R}^d:$
 \begin{itemize}
\item Euclidean distance kernel: $L_e(x,y)=\left\|x-y\right\|^r,$ where $0<r<2,$
\item Gaussian kernel: $L_g(x,y)=1-\frac{1}{2}\exp\left(-\frac{1}{2} (x-y)^{\tau} V(x-y)\right)$
\item Cauchy kernel: $L_c(x,y)=1-\frac{1}{2}\frac{1}{1+(x-y)^{\tau}  V(x-y)},
$
\end{itemize}
where $V$ stands for $d-$dimensional covariance matrix.
\end{pr}

\begin{pr}
Let $\mathcal{X}$ be a separable Hilbert space. Assume that $f(x)$ is real characteristic functional of an infinitely divisible measure on $\mathcal{X}.$ Then $L(x,y)=-\log f(x-y)$ is a negative definite kernel. We know that 
\[
\mathcal{L}(x-y)=\frac{1}{2}(B(x-y),x-y)-\int_{\mathcal{X}}\left( e^{i\left\langle x-y,t\right\rangle}-1 -\frac{i\left\langle x-y,t\right\rangle}{1+\left\| t\right\|^2}\right)\frac{1+\left\| t\right\|^2}{\left\| t\right\|^2}d\theta(t),
\]
where $B$ is the kernel operator and $\theta$ is a finite measure for which $\theta(\left\{ 0\right\})=0.$ If $\operatorname{supp} \theta =\mathcal{X},$ then $\mathcal{L}$ is strongly negative definite kernel on $\mathcal{X}.$
\end{pr}
\begin{tm}
Let $L$ be a strongly negative definite kernel on $\mathcal{X} \times \mathcal{X}$ satisfying $L(x,y)=L(y,x)$ and $L(x,x)=0$ for all $x,y \in \mathcal{X}.$
Let $\mathcal{N}: \mathcal{P}_{L} \times \mathcal{P}_{L} \to \mathbb{R}$ be defined by
$$
\mathcal{N}(\mu,\nu)=
2\int_{\mathcal X} \int_{\mathcal X}  L(x,y) d\mu (x) d\nu (y) - 
\int_{\mathcal X} \int_{\mathcal X}  L(x,y) d\mu (x) d\mu (y)
-\int_{\mathcal X} \int_{\mathcal X}  L(x,y) d\nu (x) d\nu (y).
$$
Then $\mathfrak{N}=\mathcal{N}^{1/2}$ is a distance on $\mathcal{P}_{L}.$
\label{tm:n}
\end{tm}

\section{Characteristic function of random convex compact sets}

\subsection{Space $l_0$}

Let $l_0$ denote the space of all real sequences with finite number of non-zero elements (terms). 
Since the space $l_0$ plays an important role as domain of characteristic function of random convex compact set (or more general, Fourier-Stiljes transformation of finite signed measure on $(\mathcal{K},\mathcal{B}(\mathcal{F})\cap \mathcal{K})$), we describe here a topology and the associated (Borel) measure on $l_0$.
It is used in further investigating and application of the properties of these characteristic functions. 
Less formally, $l_0$ can be seen as $\cup \mathbb{R}^n$ and it can be topologised by strict inductive limit of $\mathbb{R}^n$ with standard topology. 
An advantage of strict inductive limit topology is that it induces the same topology on $\mathbb{R}^n$ if we consider it as a subset of $l_0.$ 
The following definitions and basic properties of strict inductive limit topologies 
are taken from \cite{H89} and \cite{NB10}.

\begin{definition} 
Let $\left\{ X_n \right\}_{n \in \mathbb{N}}$ be increasing sequence of locally convex topological spaces such that each $X_n$ is a subspace of $X_{n+1}$ and $X=\cup X_n$ is a vector space.
For each $n \in \mathbb{N},$ let $\iota_n:X_n \to X$ be canonical injection of $X_n$ into $X.$ 
\textit{Strict inductive limit topology on $X$} is the finest locally convex topology in which all $\iota_n$ are continuous.
\end{definition}
Strict inductive limit topology induces the original topology on each of the component spaces $X_n$ \cite{NB10}.
\begin{prop}
Let $(X,\mathcal{T})$ be inductive limit of Hausdorff spaces $(X_n,\mathcal{T}_n)$ in which the origin has neighbourhood relatively compact in $(X_{n+1},\mathcal{T}_{n+1}).$
\begin{enumerate}
\item A set $U$ is open in $(X,\mathcal{T})$ if and only if $U\cap X_n $ is open in $(X_n,\mathcal{T}_n)$ for each $n.$
\item Let $U$ be open in $(X,\mathcal{T}).$ A map h of $U$ into topological space if continuous for the topology of $(X,\mathcal{T})$ if and only if, for each $n,$ $h|_{U \cap X_n}$ is continuous for the topology of $(X_n,\mathcal{T}_n).$ 
\end{enumerate}
\label{prop:str_prop}
\end{prop}
For the proof, see \cite{H89}.

For our purposes, we define a linear mapping (canonical injection) $\iota_{k_1,\ldots,k_n}: \mathbb{R}^n \to l_0$ by\break
$\iota_{k_1,\ldots,k_n}(x_1,\ldots,x_n)=t$ such 
as $t_{k_i}=x_i, \ i =1,\ldots,n$, and $t_j=0$ if $j \neq k_1,\ldots,k_n.$ Note that this mapping is canonical injection from $\mathbb{R}^d$ to $l_0.$
Since $l_0=\cup_{n=1}^{\infty}\iota_{1,\ldots,n}(\mathbb{R}^n)$, we can consider it to be a topological vector space with the topology of the strict inductive limit of $(\mathbb{R}^n, \mathcal{T}_n),$ where $\mathcal{T}_n$ stands for standard topology generated by 
Euclidean norm in $\mathbb{R}^n.$ Denote this topology by $\mathcal{T}.$ 
From Proposition $\ref{prop:str_prop}$, it follows that $U\in\mathcal{T}$ if and only if $\iota^{-1}_{1,\ldots,n}(U)$ is open in $\mathbb{R}^n$ for every $n\in \mathbb{N}$, and moreover $f:l_0\to \mathbb{R}$ is continuous if and only if $f\circ\iota_{1,\ldots,n}$ is continuous for every $n \in \mathbb{N}.$ 

Let $\mathcal{B}(\mathcal{T})=\sigma(\mathcal{T})$ be a Borel $\sigma$-algebra generated by $\mathcal{T}.$ In the following proposition, we provide characterization of $\mathcal{B}(\mathcal{T})-$measurable sets.

\begin{prop}
For $A\subseteq l_0$, it holds that $A \in \mathcal{B}(\mathcal{T})$ if and only if for every $n \in \mathbb{N},$  $\iota^{-1}_{1,\ldots,n}(A)\in \mathcal{B}(\mathbb{R}^n).$ 
\label{prop:sigma}
\end{prop}
\textit{Proof. }
Denote $\mathcal{G}=\left\{ A \subseteq l_0: \text{ \ for every \ } n \in \mathbb{N}, \ \iota^{-1}_{1,\ldots,n}(A)\in \mathcal{B}(\mathbb{R}^n) \right\}.$ Let us prove that $\mathcal{B}(\mathcal{T})=\mathcal{G}.$
First, we show $\mathcal{B}(\mathcal{T})\subseteq \mathcal{G}.$ 
\begin{itemize}
\item Note that $\mathcal{G}$ is $\sigma$-algebra, since
\begin{itemize}
\item $l_0\in \mathcal{G}$,
\item if $A \in \mathcal{G}$ then for every $n \in \mathbb{N},$ we have  $\iota^{-1}_{1,\ldots,n}(A)\in \mathcal{B}(\mathbb{R}^n),$ so $\left(\iota^{-1}_{1,\ldots,n}(A)\right)^c=\iota^{-1}_{1,\ldots,n}(A^c) \in \mathcal{B}(\mathbb{R}^n)$ from which it follows that $A^c \in \mathcal{G}$,
\item if $(A_k)_{k \in \mathbb{N}} \subset \mathcal{G}$ then for every $k,n \in \mathbb{N}$ $\iota^{-1}_{1,\ldots,n}(A_k)\in \mathcal{B}(\mathbb{R}^n),$ so $\cup_{k=1}^{\infty}\iota^{-1}_{1,\ldots,n}(A_k)=\iota^{-1}_{1,\ldots,n}(\cup_{k=1}^{\infty}A_k) \in \mathcal{B}(\mathbb{R}^n)$ from which it follows that $\cup_{k=1}^{\infty}A_k \in \mathcal{G}.$
\end{itemize}
\item $\mathcal{T}\subseteq \mathcal{G}$.
\end{itemize}
Therefore $\mathcal{B}(\mathcal{T})=\sigma(\mathcal{T}) \subseteq \mathcal{G}.$

Now, let us show that $\mathcal{G} \subseteq \mathcal{F}.$
Note that for $A \subset l_0$, it holds that
\[
A=A\cap l_0=\cup_{n=1}^{\infty}\left(A\cap\iota_{1,\ldots,n}(\mathbb{R}^n)\right)=\cup_{n=1}^{\infty}\iota_{1,\ldots,n}(\iota^{-1}_{1,\ldots,n}(A)).
\]
If $A \in \mathcal{G}$, then $\iota^{-1}_{1,\ldots,n}(A) \in \mathcal{B}(\mathbb{R}^n)$.
Thus $\iota_{1,\ldots,n}(\iota^{-1}_{1,\ldots,n}(A)) \in \iota_{1,\ldots,n}(\mathcal{B}(\mathbb{R}^n))=\sigma(\mathcal{T}\cap\iota_{1,\ldots,n}(\mathbb{R}^n))\subset \mathcal{B}(\mathcal{T})$.
Since $\mathcal{B}(\mathcal{T})$ is $\sigma$-algebra, it holds that $A\in \mathcal{B}(\mathcal{T}).$ 
\bigskip

Note that $\iota_{k_1,\ldots,k_n}$ are continuous and therefore it is $\mathcal{B}(\mathcal{T})-$measurable as $\iota_{k_1,\ldots,k_n}=\iota_{1,\ldots,k_n}\circ t_{k_1,\ldots,k_n},$ where  $t_{k_1,\ldots,k_n}:\mathbb{R}^{n}\to\mathbb{R}^{k_n}$ is continuous mapping defined by $t_{k_1,\ldots,k_n}(x_1,\ldots,x_n)=(y_1,\ldots,y_{k_n})$, where $y_{j}=x_i$ if $j=k_i$ and $y_{j}=0$ otherwise.

Now, we define a class of measures on $(l_0,\mathcal{B}(\mathcal{T}))$ that are used in the sequel. 
Let $\mu_g:\mathcal{B}(\mathcal{T}) \to [0,\infty]$ be defined by 
\begin{equation}
\mu_g(A)=\lim\limits_{n\rightarrow \infty}\sum\limits_{i=1}^{n}\sum\limits_{\left\{k_1,\ldots,k_i\right\} 
\subseteq \left\{1,\ldots,n\right\} }\int\limits_{\iota^{-1}_{{k_1},\ldots,{k_i}}(A)}g\circ\iota_{{k_1},\ldots,{k_i}}d\lambda_i,
\label{eq:l0_mu_g}
\end{equation}
where $\lambda_i$ is the Lebesgue measure on $\mathbb{R}^i$ and $g:l_0 \to \left[0,\infty\right\rangle $ is an arbitrary function such that 
$g\circ\iota_{{k_1},\ldots,{k_i}}$ is integrable with respect to the $i$-dimensional Lebesgue measure
 for all $i \in \mathbb{N}$ and all $\left\{k_1,\ldots,k_i\right\}\subset \mathbb{N}.$ 
Let us show that $\mu :\mathcal{B}(\mathcal{T})\to [0,\infty]$ is a measure on $(l_0,\mathcal{B}(\mathcal{T})).$ It is sufficient to show that $\mu$ is $\sigma$-additive. For an arbitrary sequence $\{A_j\}_{j \in\mathbb N}$ of disjoint sets from $\mathcal{F}$, it holds that
\begin{align*}
\mu_g(\cup_{k=1}^{\infty}A_j) &=\lim\limits_{n\rightarrow \infty}\sum\limits_{i=1}^{n}\sum\limits_{\left\{k_1,\ldots,k_i\right\} 
\subseteq \left\{1,\ldots,n\right\} }\int\limits_{\iota^{-1}_{{k_1},\ldots,{k_i}}(\cup_{k=1}^{\infty}A_j)}g\circ\iota_{{k_1},\ldots,{k_i}}d\lambda_i\\
&=\lim\limits_{n\rightarrow \infty}\sum\limits_{i=1}^{n}\sum\limits_{\left\{k_1,\ldots,k_i\right\} 
\subseteq \left\{1,\ldots,n\right\} }\sum\limits_{j=1}^{\infty}\int\limits_{\iota^{-1}_{{k_1},\ldots,{k_i}}(A_j)}g\circ\iota_{{k_1},\ldots,{k_i}}d\lambda_i=\sum\limits_{j=1}^{\infty}\mu_g(A_j).
\end{align*}

\begin{prop}
If $f:l_0 \to \mathbb{R}^+$ is an arbitrary $\mathcal{B}(\mathcal{T})-$measurable function, then
\begin{equation*}
\int\limits_{l_0}f d\mu_g=\lim\limits_{n \rightarrow \infty}\sum\limits_{i=1}^{n}\sum\limits_{\left\{k_1,\ldots,k_i\right\}
\subseteq \left\{1,\ldots,n\right\}}\int\limits_{\mathbb{R}^i}f\circ\iota_{{k_1},\ldots,{k_i}}g\circ\iota_{{k_1},\ldots,{k_i}}d\lambda_i,
\end{equation*}
where $\lambda_i$ is the $i-$dimensional Lebesgue measure.
\label{prop:mu_g}
\end{prop}
\textit{Proof. }
Let $f=\chi_{B}$ for $B \in \mathcal{B}(\mathcal{T}).$
Then 
\begin{align*}
\int fd\mu_g & =\int \chi_{B} d\mu_g=\mu_g(B)
=\lim\limits_{n\rightarrow \infty}\sum\limits_{i=1}^{n}\sum\limits_{\left\{k_1,\ldots,k_i\right\} 
\subseteq \left\{1,\ldots,n\right\} }\int\limits_{\iota^{-1}_{{k_1},\ldots,{k_i}}(B)}g\circ\iota_{{k_1},\ldots,{k_i}}d\lambda_i\\
&=\lim\limits_{n\rightarrow \infty}\sum\limits_{i=1}^{n}\sum\limits_{\left\{k_1,\ldots,k_i\right\} 
\subseteq \left\{1,\ldots,n\right\} }\int\limits_{\mathbb{R}^i} \chi_{\iota^{-1}_{{k_1},\ldots,{k_i}}(B)}g\circ\iota_{{k_1},\ldots,{k_i}}d\lambda_i\\
&=\lim\limits_{n\rightarrow \infty}\sum\limits_{i=1}^{n}\sum\limits_{\left\{k_1,\ldots,k_i\right\} 
\subseteq \left\{k_1,\ldots,k_n\right\} }\int\limits_{\mathbb{R}^i} \chi_{B}\circ \iota_{{k_1},\ldots,{k_i}}g\circ\iota_{{k_1},\ldots,{k_i}}d\lambda_i\\
&=\lim\limits_{n\rightarrow \infty}\sum\limits_{i=1}^{n}\sum\limits_{\left\{k_1,\ldots,k_i\right\} 
\subseteq \left\{1,\ldots,n\right\} }\int\limits_{\mathbb{R}^i} f\circ \iota_{{k_1},\ldots,{k_i}}g\circ\iota_{{k_1},\ldots,{k_i}}d\lambda_i\\
\end{align*}
If $f$ is simple $\mathcal{B}(\mathcal{T})-$measurable, i.e. $f=\sum_{j=1}^{J}a_j\chi_{B_j},$ where $a_1,\ldots,a_J\in \mathbb{R}$ and $B_1,\ldots,B_j\in \mathcal{B}(\mathcal{T})$ are mutually disjoint sets, then it holds that
\begin{align*}
\int f d\mu_g &=\sum_{j=1}^{J}a_j\int\chi_{B_j}d\mu_g =\sum_{j=1}^{J}a_j \left(\lim\limits_{n\rightarrow \infty}\sum\limits_{i=1}^{n}\sum\limits_{\left\{k_1,\ldots,k_i\right\} 
\subseteq \left\{1,\ldots,n\right\} }\int\limits_{\mathbb{R}^i} \chi_{B_j\circ\iota_{{k_1},\ldots,{k_i}}}g\circ\iota_{{k_1},\ldots,{k_i}}d\lambda_i\right)\\
&= \lim\limits_{n\rightarrow \infty}\sum\limits_{i=1}^{n}\sum\limits_{\left\{k_1,\ldots,k_i\right\} 
\subseteq \left\{k_1,\ldots,k_n\right\} }\int\limits_{\mathbb{R}^i} (\sum_{j=1}^{J}a_j\chi_{B_j\circ\iota_{{k_1},\ldots,{k_i}}})g\circ\iota_{{k_1},\ldots,{k_i}}d\lambda_i\\
&=\lim\limits_{n\rightarrow \infty}\sum\limits_{i=1}^{n}\sum\limits_{\left\{k_1,\ldots,k_i\right\} 
\subseteq \left\{1,\ldots,n\right\} }\int\limits_{\mathbb{R}^i} f\circ \iota_{{k_1},\ldots,{k_i}}g\circ\iota_{{k_1},\ldots,{k_i}}d\lambda_i.
\end{align*}
If $f$ is an arbitrary non negative $\mathcal{B}(\mathcal{T})-$measurable function, then there exists a non decreasing sequence $\left\{f_j\right\}_{j\in \mathbb{N}}$ of simple functions such that $f(t)=\lim\limits_{
j\rightarrow \infty}f_j(t),$ $t \in l_0.$ Then from the Lebesgue theorem of monotone convergence, we get 
\begin{align*}
\int & f d\mu =\lim\limits_{j \rightarrow \infty}f_j d\mu=\lim\limits_{j \rightarrow \infty}\lim\limits_{n\rightarrow \infty}\sum\limits_{i=1}^{n}\sum\limits_{\left\{k_1,\ldots,k_i\right\} 
\subseteq \left\{1,\ldots,n\right\} }\int\limits_{\mathbb{R}^i} f_j\circ \iota_{{k_1},\ldots,{k_i}}g\circ\iota_{{k_1},\ldots,{k_i}}d\lambda_i\\
&=\lim\limits_{n\rightarrow \infty}\sum\limits_{i=1}^{n}\sum\limits_{\left\{k_1,\ldots,k_i\right\} 
\subseteq \left\{1,\ldots,n\right\} }\int\limits_{\mathbb{R}^i} \lim\limits_{j \rightarrow \infty}f_j\circ \iota_{{k_1},\ldots,{k_i}}g\circ\iota_{{k_1},\ldots,{k_i}}d\lambda_i\\
&=\lim\limits_{n\rightarrow \infty}\sum\limits_{i=1}^{n}\sum\limits_{\left\{k_1,\ldots,k_i\right\} 
\subseteq \left\{1,\ldots,n\right\} }\int\limits_{\mathbb{R}^i} f\circ \iota_{{k_1},\ldots,{k_i}}g\circ\iota_{{k_1},\ldots,{k_i}}d\lambda_i.
\end{align*}
\bigskip

We conclude this section by two examples of $\sigma$-finite measure  on $(l_0,\mathcal{B}(\mathcal{T}))$ of the form defined by the relation $(\ref{eq:l0_mu_g}).$ 
\begin{pr}
If  $A_n=\iota_{1,\ldots,n}(\mathbb{R}^n)$, i.e. $A_n$ is set of all sequences that have only non-zero elements at the first $n$ positions, then $\left\{A_n\right\}$ is increasing and $\cup A_n=l_0.$
Let $$g(t)=\operatorname{exp}\left(-\sum\limits_{n=1}^{\infty}\left|t_i\right|/w_i\right),$$ 
where $\left\{w_i\right\}_{i \in \mathbb{N}}$ is an arbitrary sequence of non negative real numbers. It is easy to see that
$$\mu_g(A_n)=\sum\limits_{i=1}^{n}\sum\limits_{\left\{i_1,\ldots,i_k\right\}\subseteq \left\{1,\ldots,n\right\}}w_{i_1}\cdots w_{i_k}<\infty$$
and
\[
\mu_g(l_0)=\lim\limits_{n\rightarrow \infty}\mu_g(A_n)=\infty.
\]
Thus the measure $\mu_g$ is $\sigma$-finite and $\operatorname{supp}\mu_g=l_0,$ since $\mu_g(U)>0$ for every $U \in \mathcal{T}.$
\label{ex:g}
\end{pr}

\begin{pr}
Denote $c_r=\pi^{\frac{r+1}{2}}/ \Gamma\left(\frac{r+1}{2}\right).$ For $0<r<2$, let
$$h(t)=c_r\left( \sum\limits_{n=1}^{\infty}t_n^2 \right)^{-\frac{r+1}{2}}.$$ 
Define $A_{n,m}=\iota_{1,\ldots,n}\left(B^n(0,m)\diagdown B^n(0,1/m)\right),$ where $B^n(0,m)$ denotes the open ball in $\mathbb{R}^n$ with the centre in the origin and the radius $m.$  It holds that $l_0=\cup_{n\in \mathbb{N}}\cup_{m\in \mathbb{N}}A_{n,m}\cup \left\{ 0 \right\}$ and $\mu_h(A_{n,m})<\infty,$ so again, the measure $\mu_h$ is $\sigma$-finite and $\operatorname{supp}\mu_h=l_0$ as $\mu_h(U)>0$ for every $U \in \mathcal{T}.$
\label{ex:h}
\end{pr}

\subsection{Fourier-Stieltjes transform of finite signed measure on $(\mathcal{K},\mathcal{U})$}

This section introduces the Fourier-Stieltjes transform of a finite signed measure on $(\mathcal{K},\mathcal{U}).$ 

\begin{definition} 
Let $\nu$ be a finite signed measure on $(\mathcal{K},\mathcal{U})$ and $\left\{u_j: j \in \mathbb{N}\right\}$ be a countable dense subset of $S^{d-1}.$ 
The \textit{Fourier-Stieltjes transform of $\nu$} is the function
$\hat{\nu}: l_0 \to \mathbb{C}$ defined by
\[
\hat{\nu}(t)=\int\limits_{\mathcal{K}}e^{i\sum\limits_{j=1}^{\infty}t_j h(A,u_j)} d\nu(A).
\]
If $\nu$ is a probability measure, then its Fourier-Stieltjes transformation is called \textit{characteristic function of the measure $\nu.$}
\end{definition}

\begin{tm}[Uniqueness theorem]
Let $\nu_1$ and $\nu_2$ be a finite signed measures on $(\mathcal{K},\mathcal{U}).$ It holds that
\[
\nu_1=\nu_2 \Longleftrightarrow \hat{\nu_1}=\hat{\nu_2}.
\]
\label{tm:uniqueH}
\end{tm}
\textit{Proof. }
It is sufficient to show that $\hat{\nu}=0$ if and only if $\nu=0.$
If $\nu=0$, then it is obvious that $\hat{\nu}=0.$ 
Let us show the converse.
For arbitrary set of indices $I=\left\{j_1,\ldots,j_n\right\} \subset \mathbb{N}$, let
\[
\nu_I(B)=\nu(\left\{A \in \mathcal{K}: \left(h(A,u_{j_1}),\ldots,h(A,u_{j_n})\right) \in B\right\}), \ B \in \mathcal{B}(\mathbb{R}^n).\]
It is obvious that $\nu_I$ is finite signed measure on $(\mathbb{R}^n,\mathcal{B}(\mathbb{R}^n)).$\\
For $t \in \mathbb{R}^n$, let $t'=\iota_{j_1,\ldots,j_n}(t) \in l_0.$
Since $\hat{\nu}_I(t)=\hat{\nu}(t')=0$ holds, it follows from Theorem $\ref{tm:unique}$ that $\nu_I=0.$
Since the algebra of the sets $\left\{A \in \mathcal{K}: \left(h(A,u_{j_1}),\ldots,h(A,u_{j_n})\right) \in B\right\}$ generates $\sigma$-algebra $\mathcal{U}$ and $\nu(\left\{A \in \mathcal{K}: \left(h(A,u_{j_1}),\ldots,h(A,u_{j_n})\right) \in B\right\})=0$, we get $\nu=0.$
\bigskip

\section{Test of equality in distribution of random convex compact sets}
\label{sec:Test}
Consider two groups of random convex compact sets, where each group is coming from the same distribution. In this section, we introduce two approaches for testing the equality of the distributions of these two groups using the theory from the previous sections.

The tests are based on $\mathfrak{N}$-distances.
First, in Section $\ref{sec:Ndist}$, we introduce theoretical results used in the construction of the tests, more precisely, we construct the class of $\mathfrak{N}$-distances on the space of distributions of random convex compact sets in order to create the test statistics. 
In Section $\ref{sec:L}$, we deal with evaluation of  $\mathfrak{N}$-distance part based on available information from the sample.
Finally in Section $\ref{sec:test}$, we present two general procedures for formulating distribution free two-sample tests based on $\mathfrak{N}$-distances (see\cite{K06}). 

\subsection{$\mathfrak{N}$-distance on distributions of random convex compact sets}
\label{sec:Ndist}

Let $\left\{u_n\right\}_{n \in \mathbb{N}}$ be a countable dense subset of the unit sphere $S^{d-1} \subset \mathbb{R}^{d}.$
For $A \in \mathcal{K}$ and $t \in l_0$, we define  $\varphi_t :  \mathcal{K} \to \mathbb{C}$ for arbitrary $t \in l_0$ as
$$
\varphi_t(A)=\operatorname{exp}\left(i\sum\limits_{n=1}^{\infty}t_n h(A,u_n)\right).
$$ 
The values of $\varphi_t(A)$ are always finite since the number of non-zero terms in the sum is finite so the value of the sum is finite.
For $t \in l_0$, we define $K_t:\mathcal{K}\times \mathcal{K} \to \mathbb{C}$ by
\begin{equation}
K_t(A,B)=\frac{1}{2}\varphi_t(A)\overline{\varphi_t(B)},
\end{equation}
Proposition $\ref{prop:pdk}$ shows that $K_t$ is positive definite kernel.
From Proposition $\ref{prop:ndk}$, it follows that the function $L_t:\mathcal{K}\times \mathcal{K} \to \mathbb{C}$ defined by
\begin{align*}
L_t(A,B)&=K_t(A,A)+K_t(B,B)-K_t(A,B)-K_t(B,A)=\frac{1}{2}\left|\varphi_t(A)-\varphi_t(B)\right|^2\\
&=\frac{1}{2}\left|\operatorname{exp}\left(i\sum\limits_{n=1}^{\infty}t_n h(B,u_n)\right)\right|^2\left|1-\operatorname{exp}\left(i\sum\limits_{n=1}^{\infty}t_n (h(A,u_n)-h(B,u_n))\right)\right|^2\\
&=\frac{1}{2}\left(1-\cos\left(\sum\limits_{n=1}^{\infty}t_n \left(h(A,u_n)-h(B,u_n)\right)\right)\right)^2+\frac{1}{2}\sin^2\left(\sum\limits_{n=1}^{\infty}t_n (h(A,u_n)-h(B,u_n))\right)\\
&=1-\cos\left(\sum\limits_{n=1}^{\infty}t_n (h(A,u_n)-h(B,u_n))\right).
\end{align*}
is a negative definite kernel on $\mathcal{K}.$
Obviously, it holds that
\[
L_t(A,B)=0  \ \forall t \in l_0 \ \Leftrightarrow \ A=B.
\]

The following theorem introduces the generalised version of $\mathfrak{N}$-distance that is allowed to take infinite values and whose domain is the whole $\mathcal{P}$ .
\begin{tm}
\label{tm:ndist}
Let $\mu$ be a $\sigma$-finite measure on $l_0$ with $\operatorname{supp}\mu=l_0$ and $\left( A_n\right)_{n \in \mathbb{N}}$ a sequence of sets in $l_0$ such that $l_0=\cup_{n\in \mathbb{N}}A_n,$ $A_n\subset A_{n+1}$ and $\mu(A_n)<\infty$ for all $n \in \mathbb{N}.$ Define function $\mathcal{N}: \mathcal{P} \times \mathcal{P} \to \overline{\mathbb{R}}$ by
\begin{align*}
\mathcal N(\nu_1,\nu_2) = & \lim\limits_{n \rightarrow \infty}\left(  2\int_{\mathcal{K}} \int_{\mathcal{K}} \int_{A_n}L_t(A,B)d\mu(t) d\nu_1 (x) d\nu_2 (y)
- \int_{\mathcal{K}} \int_{\mathcal{K}} \int_{A_n}L_t(A,B)d\mu(t) d\nu_1 (x) d\nu_1 (y) \right.\\
& \left.
-\int_{\mathcal{K}} \int_{\mathcal{K}} \int_{A_n}L_t(A,B)d\mu(t) d\nu_2 (x) d\nu_2 (y) \right).
\end{align*}
Then $\mathfrak{N}=\mathcal N^{1/2}$ is a metric on the space $\mathcal{P}$ of all distributions of random convex compact sets.
\end{tm}
\textit{Proof. }
Since the measure $\mu$ is $\sigma$-finite, there exists a sequence of the sets $A_n\subset l_0$ such that $l_0=\cup_{n\in \mathbb{N}}A_n,$ $A_n\subset A_{n+1}$ and $\mu(A_n)<\infty$ for all $n \in \mathbb{N}.$
Let us define
\begin{equation}
\label{eq:Ln}
\mathcal{L}_n(A,B)=\int_{A_n}L_t(A,B)d\mu(t).
\end{equation}
It holds that $\mathcal{L}(A,B)=\lim\limits_{n \rightarrow \infty}\mathcal{L}_n(A,B).$
Note that  $0 \leq L_t(A,B) \leq 4$ for all $t \in l_0.$ It implies that $0 \leq \mathcal{L}_n(A,B) \leq 4\mu(A_n)<\infty$ for all $n \in\mathbb{N}$ and so $\mathcal{P}_{L_n}=\mathcal{P}$ for all $n \in \mathbb{N}.$

Further, we define $\mathcal{N}_n :\mathcal P \times \mathcal P \to \mathbb{R}$ by
\begin{align}
\label{eq:Nn}
\mathcal N_n(\nu_1,\nu_2) &= 2\int_{\mathcal{K}} \int_{\mathcal{K}} \mathcal{L}_n(A,B) d\nu_1 (A) d\nu_2 (B) - 
\int_{\mathcal{K}} \int_{\mathcal{K}} \mathcal{L}_n(A,B) d\nu_1 (A) d\nu_1 (B)
\\
&
-\int_{\mathcal{K}} \int_{\mathcal{K}} \mathcal {L}_n(A,B) d\nu_2 (A) d\nu_2 (B). \nonumber
\end{align}
 Since the difference $\mathcal {L}_{n+1}(x,y)-\mathcal {L}_n(x,y)$ is a negative definite kernel as an integral of negative definite kernels,  it follows from Proposition $\ref{prop:strong}$ that $\mathcal{N}_n(\nu_1,\nu_2)\leq\mathcal{N}_{n+1}(\nu_1,\nu_2).$ It is easy to see that
\begin{equation}
\mathcal{N}(\nu_1,\nu_2)=\lim\limits_{n \rightarrow \infty}\mathcal{N}_n(\nu_1,\nu_2)
\label{eq:Nlim}
\end{equation}
and it follows that $\mathcal{N}$ is well defined as a limit of non decreasing sequence. It holds that $$\mathcal{N}_n(\nu_1,\nu_2)=\int_{A_n}\left|\hat{\nu_1}-\hat{\nu_2}\right|^2 d\mu,$$
so 
$$\mathcal{N}(\nu_1,\nu_2)=\int_{l_0}\left|\hat{\nu_1}-\hat{\nu_2}\right|^2 d\mu,$$
and from this relation, it follows that  $\mathfrak{N}=\mathcal{N}^{1/2}$ is a distance on the set $\mathcal{P}$ of all probability measures on $(\mathcal{K},\mathcal{U}).$
\bigskip

Let $\mathbf{A}$ and $\mathbf{B}$ be random convex compact sets with the distributions $\nu^{\mathbf{A}}$ and $\nu^{\mathbf B},$ respectively. 
Denote $\hat{\nu}_{k_1,\ldots,k_j}^{\mathbf A}$ and $\hat{\nu}_{k_1,\ldots,k_j}^{\mathbf B}$ characteristic functions of random vectors\break $(h( \mathbf A,u_{k_1}),\ldots,h(\mathbf A,u_{k_j}))$ and $(h(\mathbf B,u_{k_1}),\ldots,h(\mathbf B,u_{k_j}))$, respectively.
If we consider the measure $\mu$ to be 
$\mu_g$ defined by the relation $(\ref{eq:l0_mu_g})$ and $A_n=\iota_{1,\ldots,n}(\mathbb{R}^n)$, it is obvious that 
\begin{equation}
\mathcal{N}_n(\nu^{\mathbf A},\nu^{\mathbf B})=\sum\limits_{j=1}^{n}\sum\limits_{\left\{k_1,\ldots,k_j\right\}
\subseteq \left\{1,\ldots,n\right\}}\int\limits_{\mathbb{R}^j}\left|\hat{\nu}_{k_1,\ldots,k_j}^{\mathbf A}-\hat{\nu}_{k_1,\ldots,k_j}^{\mathbf B}\right|^2g\circ\iota_{{k_1},\ldots,{k_j}}d\lambda_j.
\label{eq:Ndist_char}
\end{equation}
 So $\mathcal{N}_n(\nu^{\mathbf A},\nu^{\mathbf B})=0$ implies that all finite dimensional distributions of random vectors 
$\left(h_{\mathbf A}(u_1),\ldots,h_{\mathbf A}(u_n)\right)$ and $\left(h_{\mathbf B}(u_1),\ldots,h_{\mathbf B}(u_n)\right)$ are equal, and $\mathcal{N}(\nu^{\mathbf A},\nu^{\mathbf B})=0$ implies equality in distribution of the
support functions seen as two random processes on $S^{d-1}$.

\subsection{Approximation of $\mathfrak{N}$-distance based on sample}
\label{sec:L}

Here we give a proposal how to calculate truncated version of $\mathfrak{N}$-distance introduced in Theorem \ref{tm:ndist} using the support functions from realisations of convex compact sets evaluated at finite set of directions $u_1,\ldots,u_n \in S^{d-1},$ $ n \in \mathbb{N}.$

From the relations $(\ref{eq:Nn})$ and $(\ref{eq:Nlim})$, it follows that the mentioned $\mathfrak{N}$-distance can be expressed as a limit of integrals of negative definite kernel $\mathcal{L}_n,$ $n \in \mathbb{N}$ defined by $(\ref{eq:Ln}).$
Suppose that the measure $\mu$ on $l_0$ takes the form of $\mu_g$ defined by relation $(\ref{eq:l0_mu_g})$ and that  $A_n=\iota_{1,\ldots,n}(\mathbb{R}^n),$ $n \in \mathbb{N}.$ 
Let $A$ and $B$ be two realisations of random convex compact sets, and suppose we know the values of their support functions $h(A,u_1),\ldots,h(A,u_n)$ and $h(B,u_1),\ldots,h(B,u_n)$, respectively.
From Proposition $\ref{prop:mu_g}$, we have that $\mathcal{L}_n$ is of the form
\begin{equation}
\mathcal{L}_n(A,B)=\sum\limits_{j=1}^{n}\sum\limits_{\left\{k_1,\ldots,k_j\right\} 
\subseteq \left\{1,\ldots,n\right\} }\int\limits_{\mathbb{R}^j}\left(1-\cos\left(\sum\limits_{l=1}^{j}a_{k_l}t_{l}\right)\right)g\circ\iota_{k_1,\ldots,k_j} dt_1\ldots dt_j,
\label{eq:l}
\end{equation}
where $a_j= h(A,u_{j})-h(B,u_j)$ for $j=1,\ldots,n.$
For $n \in \mathbb{N}$ and $1\leq D\leq n,$ we define $\mathcal{L}_{(n,D)}(A,B)$ as further approximation of $\mathcal{L}_n(A,B)$  from $(\ref{eq:l})$ as
\begin{equation}
\mathcal{L}_{(n,D)}(A,B)=\sum\limits_{j=1}^{D}\sum\limits_{\left\{k_1,\ldots,k_j\right\} 
\subseteq \left\{1,\ldots,n\right\} }\int\limits_{\mathbb{R}^j}\left(1-\cos\left(\sum\limits_{l=1}^{j}a_{k_l}t_{l}\right)\right)g\circ\iota_{k_1,\ldots,k_j} dt_1\ldots dt_j,
\label{eq:lhat_D}
\end{equation}
where $a_j= h(A,u_{j})-h(B,u_j)$ for $j=1,\ldots,n.$
Now, we define 
\begin{align}
\label{eq:N_nD}
\mathfrak N_{(n,D)}(\nu_1,\nu_2) &= 2\int_{\mathcal{K}} \int_{\mathcal{K}} \mathcal{L}_n(A,B) d\nu_1 (A) d\nu_2 (B) - 
\int_{\mathcal{K}} \int_{\mathcal{K}} \mathcal{L}_n(A,B) d\nu_1 (A) d\nu_1 (B)
\\
&
-\int_{\mathcal{K}} \int_{\mathcal{K}} \mathcal {L}_n(A,B) d\nu_2 (A) d\nu_2 (B) \nonumber
\end{align}
as a approximation of $\mathfrak{N}$-distance based on information available.

The value of $D$ was introduced in order to reduce the difficulty and the time-consumption of calculations, and moreover, it allows to choose the depth of investigating the dependence structure of $A$ and $B$ as  $(\ref{eq:Ndist_char})$ implies equality of finite dimensional distributions for dimensions less or equal to $D$. 
Note that $D=n$ corresponds to testing equality of distributions of random vectors $(h(A,u_1),\ldots,h(A,u_n))$ 
and  $(h(B,u_1),\ldots,h(B,u_n))$.

There is a lot of possibilities how to choose the weight functions $g:l_0 \to \mathbb{R}$.
In the following two examples, we present evaluation of $\mathcal{L}_{(n,D)}$ from $(\ref{eq:l_D})$, where the functions $g$'s are taken from the examples $\ref{ex:g}$ and $\ref{ex:h}.$
\begin{pr}
\label{ex:new_kernel1}
Consider the function $g$ from Example $\ref{ex:g}.$ The following lemma helps us to calculate the integrals in $(\ref{eq:l})$ by recursion. 
\begin{lm}
Let $a_1,\ldots,a_n$ be an arbitrary sequence of real numbers and $w_1,\ldots,w_n$ be an arbitrary sequence of positive real numbers.
Define
\[
I_c(k)=\int\limits_0^{\infty}\ldots\int\limits_0^{\infty}\cos\left(\sum\limits_{j=1}^{k}a_jt_j\right)\exp\left(-\sum\limits_{j=1}^{k}(t_j/w_j)\right)dt_1\ldots dt_k\]
and
\[
I_s(k)=\int\limits_0^{\infty}\ldots\int\limits_0^{\infty}\sin\left(\sum\limits_{j=1}^{k}a_jt_j\right)\exp\left(-\sum\limits_{j=1}^{k}(t_j/w_j)\right)dt_1\ldots dt_k.
\]
It holds that
$$
I_c(1)
=\frac{1}{1+(w_1a_1)^2}, 
\quad \quad
I_s(1)
=\frac{w_1^2a_1}{1+(w_1a_1)^2}
$$
and
\[
\left[\begin{matrix}
I_c(k)\\
I_s(k)
\end{matrix}\right] =(w_k/(1+(w_ka_k)^2))\left[ \begin{matrix} 1 & -w_k a_k\\ w_k a_k & 1 \end{matrix}\right] \left[\begin{matrix}
I_c(k-1)\\
I_s(k-1)
\end{matrix}\right] .
\]
\end{lm}
\textit{Proof. }
This Lemma is direct consequence of the fact that
\begin{align*}
\int\limits_0^{\infty}\cos(b+at)\operatorname{exp}(-t/w)dt&=w(\cos(b)-wa\sin(b))/(1+w^2a^2),\\
\int\limits_0^{\infty}\sin(b+at)\operatorname{exp}(-t/w)dt&=w(aw\cos(b)+\sin(b))/(1+w^2a^2).
\end{align*}
\bigskip
\end{pr}
\begin{pr}
\label{ex:new_kernel2}
Consider the function $h$ from Example $\ref{ex:h}.$ Then for $0< r <2,$ 
\begin{equation}
\mathcal{L}_{(n,D)}(A,B)=\sum\limits_{j=1}^{D}\sum\limits_{\left\{k_1,\ldots,k_j\right\} 
\subseteq \left\{1,\ldots,n\right\} }\left(\sum\limits_{l=1}^j \left(h(A,u_{k_l})-h(B,u_{k_l})\right)^2 \right)^{r/2},
\label{eq:l_D}
\end{equation}
as $\left\|x\right\|=c_r\int\left( 1-\cos\left( \left\langle t,x\right\rangle\right)\right)/\left\|t\right\|^{r+1}dt,$ where $x,t \in\mathbb{R}^n$ and $\left\| \cdot \right\|$ is the Euclidean norm.

\end{pr}
\subsection{Empirical estimates and test of equality in distribution}
\label{sec:test}

In this section, we show minor modifications of two standard approaches for testing equality in distribution of random elements using $\mathfrak{N}$-distances (see\cite{K06}) for our case of random convex compact sets.

Suppose we have $m_1$ support functions $h_1^{(1)},\ldots,h_{m_1}^{(1)}$ of realisations $\mathbf A_1,\ldots, \mathbf A_{m_1}$ from the random convex compact set $ \mathbf A$ with the distribution $\nu_1$
and $m_2$ support functions $h_1^{(2)},\ldots,h_{m_2}^{(2)}$ of realisations $\mathbf B_1,\ldots, \mathbf B_{m_2}$ from the random convex compact set $\mathbf B$ with the distribution $\nu_2$, all of them evaluated in $n$ 
different directions $u_1,\ldots,u_n \in S^{d-1}.$
We want to test whether the samples come from the same distribution without any assumption on their distributions, so we are addressing so called two-sample problem. 
We test more general null hypothesis $H_0:$ finite dimensional distributions of random vectors $(h(\mathbf{A},u_1),\ldots,h(\mathbf{A},u_n))$ 
and  $(h(\mathbf{B},u_1),\ldots,h(\mathbf{B},u_n))$ are equal for dimensions less or equal to $D$, $1\leq D\leq n,$ against the alternative hypothesis $H_A:$ the distributions are not the same.

Let $\mathcal{L}_{(n,D)}$ be a negative definite kernel on $\mathcal{K} \times \mathcal{K}$ defined by $(\ref{eq:lhat_D})$ and $\mathfrak{N}_{(n,D)}$ is introduced in $(\ref{eq:N_nD}).$ 
When the samples are of the same size $m=m_1=m_2,$ the first approach is to split each of those two samples randomly into three sub-samples $\mathbf A',\mathbf A'', \mathbf A'''$ and $\mathbf B',\mathbf B'',\mathbf B'''$, respectively, of the size $m/3$ (assuming it is an integer) and define two independent samples of random variables $U_1,\ldots,U_{m/3}$ and $V_1,\ldots,V_{m/3}$ by
\begin{align*}
U_i &= \mathcal{L}_{(n,D)}(\mathbf{A'_i},\mathbf{B'_i}) - \mathcal{L}_{(n,D)}(\mathbf{A'_i},\mathbf{A''_i}),\\
V_i &= \mathcal{L}_{(n,D)}(\mathbf{B''_i},\mathbf{B'''_i}) - \mathcal L_{(n,D)}(\mathbf{A'''_i},\mathbf{B'''_i}),
\end{align*}
$i=1,\ldots,m/3$.
Under the null hypothesis these two samples of random variables are equally distributed while under the alternative hypothesis they are not equally distributed since $\mathbb{E}U-\mathbb{E}V= \mathfrak{N}_{n,D} (\mu,\nu)\neq 0.$
Now, a common univariate two sample tests can be used, for example the Kolmogorov-Smirnov test (see \cite{M51}) or the Anderson-Darling test (see \cite{P76}), denoted as KS and AD in the sequel. 

For small values of $m$ ($m<120),$ the loss of information caused by splitting the files to the size $m/3$ leads to the situation where the use of asymptotic statistics for the Kolmogorov-Smirnov test is not recommended. In the case of smaller $m$ or in the case of $m_1 \neq m_2$, it is better to use permutation version of the test. It means that first, we evaluate the empirical unbiased U-statistics estimate (see \cite{S80}) of the square of $\mathfrak{N}_{(n,D)}$
\begin{align*}
\hat{\mathcal{N}}_{(n,D)}=&\frac{2}{m_1 m_2}\sum\limits_{i=1}^{m_1}\sum\limits_{j=1}^{m_2}\mathcal{L}_{(n,D)}(\mathbf A_i,\mathbf B_j)
-\frac{1}{m_1 (m_1-1)}\sum\limits_{i=1}^{m_1}\sum\limits_{j=1}^{m_1}\mathcal{L}_{(n,D)}(\mathbf A_i,\mathbf A_j) 
\\
&
-\frac{1}{m_2(m_2-1)}\sum\limits_{i=1}^{m_2}\sum\limits_{j=1}^{m_2}\mathcal{L}_{(n,D)}(\mathbf B_i,\mathbf B_j).
\end{align*}
Then we choose the number of permutations $s$ (about 1000 is recommended), and $s$ times, we permute the considered
$1,\ldots,m_1+m_2$ sets. Consequently we split them back into two samples of the length $m_1$ and $m_2$, respectively, and for each of the permutations, we calculate $\hat{\mathcal{N}}^{(i)}_{(n,D)},$ $i=1,\ldots,s.$
Under the null hypothesis, the permutations do not modify the distribution of random variable $\mathfrak{N}.$ 
If the distributions of the two samples differ, we expect that after the permutations, the value of $\mathfrak{N}$-distance is smaller. Thus, we define
\[
p=\frac{\#\left\{i:  \hat{\mathcal{N}}^{(i)}_{(n,D)}\geq\hat{\mathcal{N}}_{(n,D)}\right\}+1}{s+1}
\]
which is the estimate of the probability that under the null hypothesis, the $\mathfrak{N}$-distance is larger than its measured value.
For a given significance level $\alpha \in \left[0,1\right],$ we reject the null hypothesis if $p\leq \alpha$, otherwise the null hypothesis is not rejected.

Note that in both proposed test, kernel $\mathcal{L}_{(n,D)}$ could be replaced by an arbitrary strongly negative definite kernel on $\mathbb{R}^{n},$ e.g. the ones from Example \ref{ex:R_d_kernels}, which would lead to standard kernel test for equality in distribution of random vectors.

\section{Assessing similarity of general random sets}

In this section, we show how to use the test for equality of distribution of the convex compact sets for assessing similarity of realisations of general random sets.

\subsection{Two-realisation problem}
\label{sec:tworeal}

Suppose that we have two realisations of random sets. 
The aim is to decide whether they come from the same underlying process. It is based on comparing their inner structure, which can be considered as a union of convex compact sets.
The approach is as follows.
\begin{enumerate}
\item Approximate each of the realisations by a union of convex compact sets with disjoint interiors in the following way: 
\begin{enumerate}
\item cover each of the realisations by discs with the same radii so that the disc centres form a Poisson disc sample process (see \cite{Eetal11}),
\item construct the Voronoi tessellation over the disc centres (see \cite{Cetal13}),
\item the cells of the Voronoi tessellation intersected with the covering discs form the division of the realisation approximation by unions of convex compact sets.
\end{enumerate}
\item From each of the tessellations over the disc union, randomly sample a given number of cells.
(note that they are convex compact).
\item Apply the test described in Section \ref{sec:test}. 
The obtained $p$-value is the measure of similarity of the realisations.
\end{enumerate}

\begin{figure}[t]
\begin{center}
\includegraphics[width=15cm]{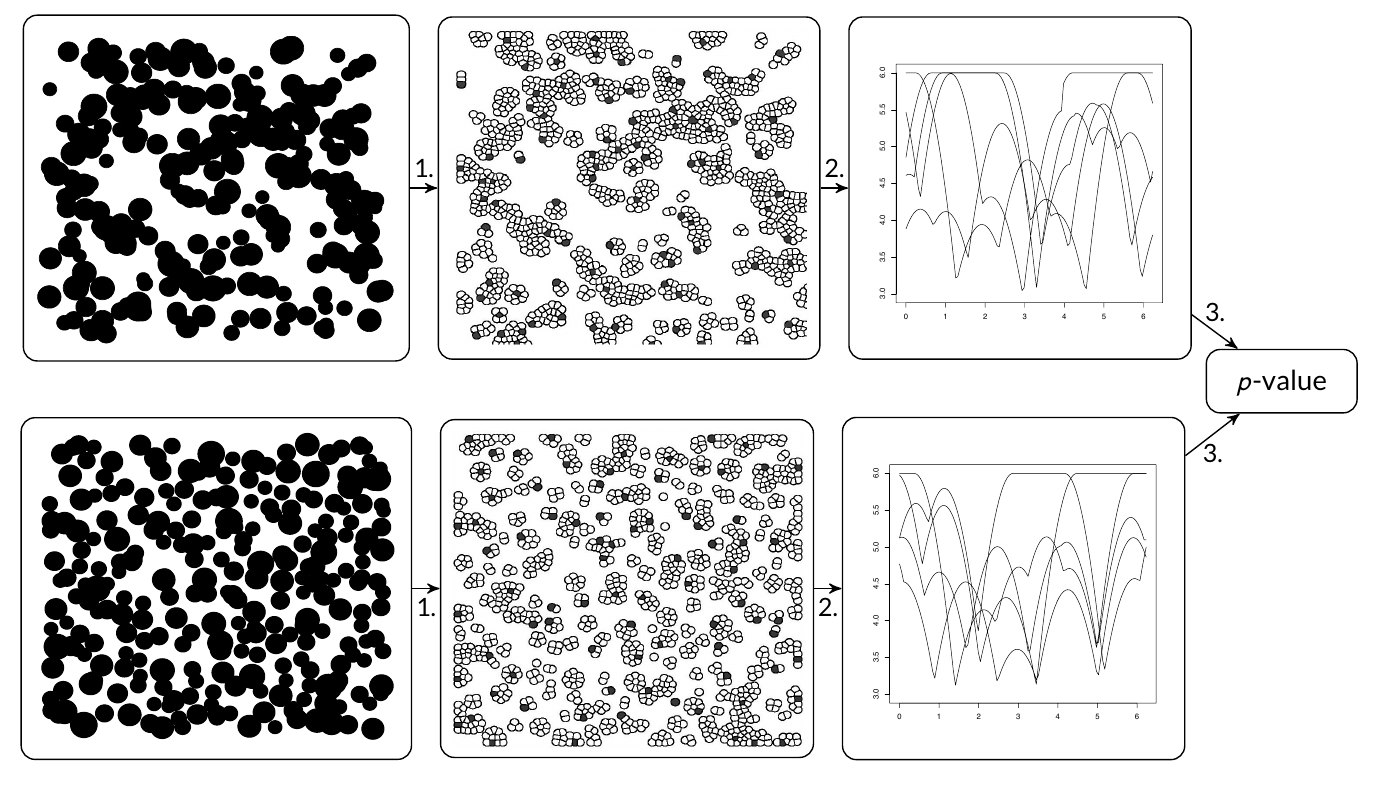}
\end{center}
\caption{
Steps of testing equality in distribution of random sets when having two realisations only.
}
\label{fig:steps}
\end{figure}

For better imagination, the procedure is graphically shown in Figure \ref{fig:steps}.
Note that the radius of covering discs in step 1. must be the same for both realisations. It is chosen so that the approximation of both the realisations is precise enough, i.e. we take the smaller of the optimal radii derived for the realisations.
More details concerning the approximation can be found in \cite{Getal16}. 

\subsection{Simulation study}

\begin{figure}[t]
\begin{center}
\includegraphics[width=15cm]{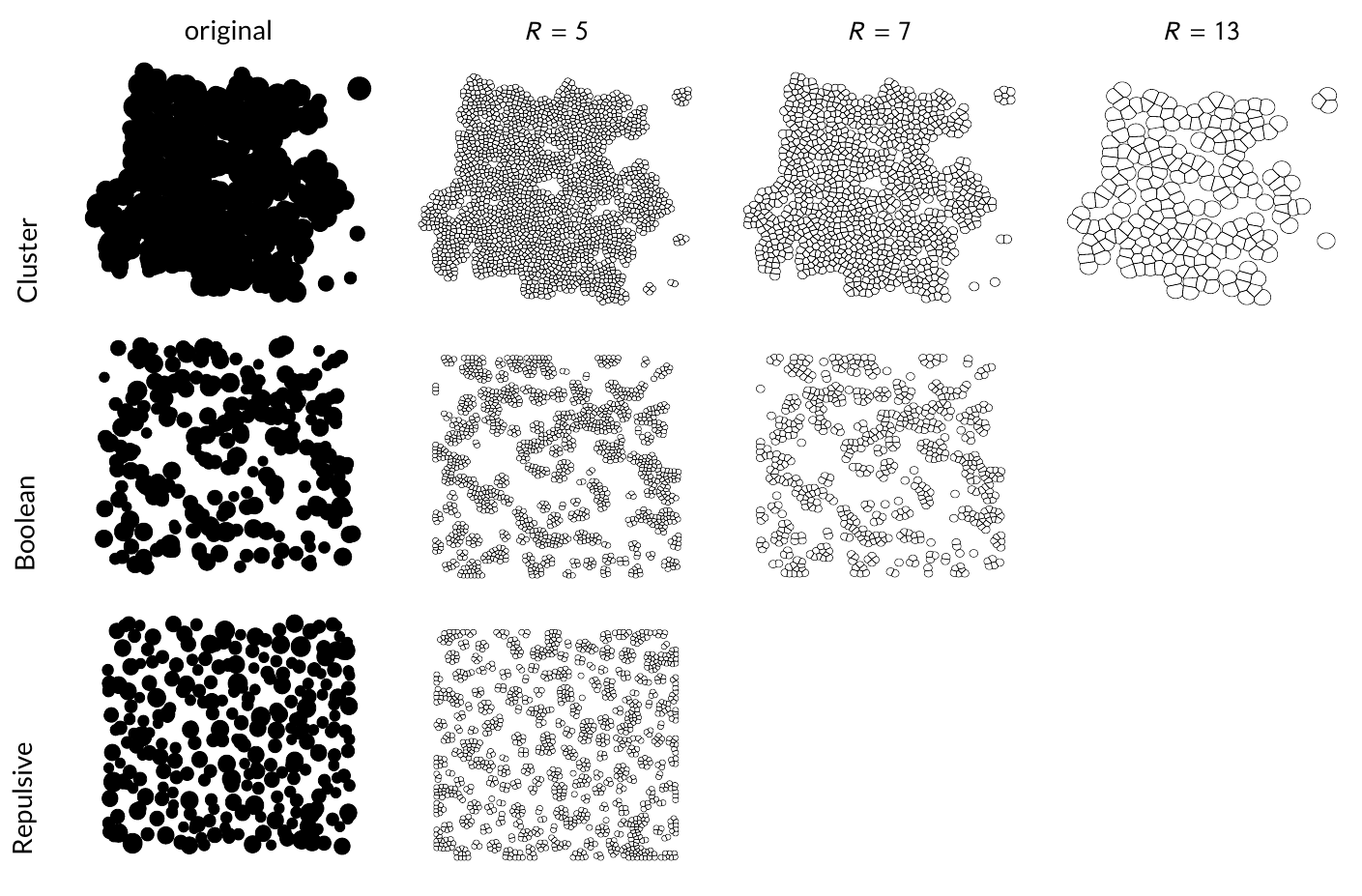}
\end{center}
\caption{Examples of realisations of the cluster process, the Boolean model and the repulsive process,
and their approximations 
used for testing.}
\label{tab:radii2}
\end{figure}

In our simulation study, we consider three different random set models.
The first one is a random-disc Boolean model (see \cite{Cetal13}), where the centres of discs form a Poisson point process and discs radii are independent on each other as well as on the disc centres. 
The second process is a disc process whose discs form clusters, therefore we call it a cluster process in the sequel. 
The third process is called repulsive process since the discs tend to be separated from each other. 
The realisations of the cluster and repulsive process are obtained by the simulations of Quermass-interaction process with suitably chosen parameters (see \cite{MH08}). 
Examples of the realisations of these three processes together with their approximations described in Section~\ref{sec:tworeal} are shown in Figure~\ref{tab:radii2}.
We consider 400 realisations of the Boolean model, cluster process and repulsive process, respectively,
and apply the test using 200 and 200 realisations  of the same processes 
as well as 200 and 200 realisations of different processes.
Thus, we have six different combinations of 200 pairs in the study.

We plot the original images in resolution 400$\times$400 pixels, then the optimal radii are $R=5$ pixels for the repulsive process, $R=7$ pixels for the Boolean model and $R=13$ pixels for the cluster process. 
It means that when we compare two realisations of the cluster processes, we use the radius $R=13$, when we compare a realisation of the cluster processes and a realisation of the Boolean model, we use $R=7$, and when we compare realisations of the cluster process and the repulsive process, we use $R=5$.
Analogously, we use $R=7$ for comparing two Boolean models and $R=5$ for comparing the Boolean model and the repulsive process as well as for comparing two repulsive processes.

\begin{figure}[t]
    \centering
\includegraphics[width=15cm]{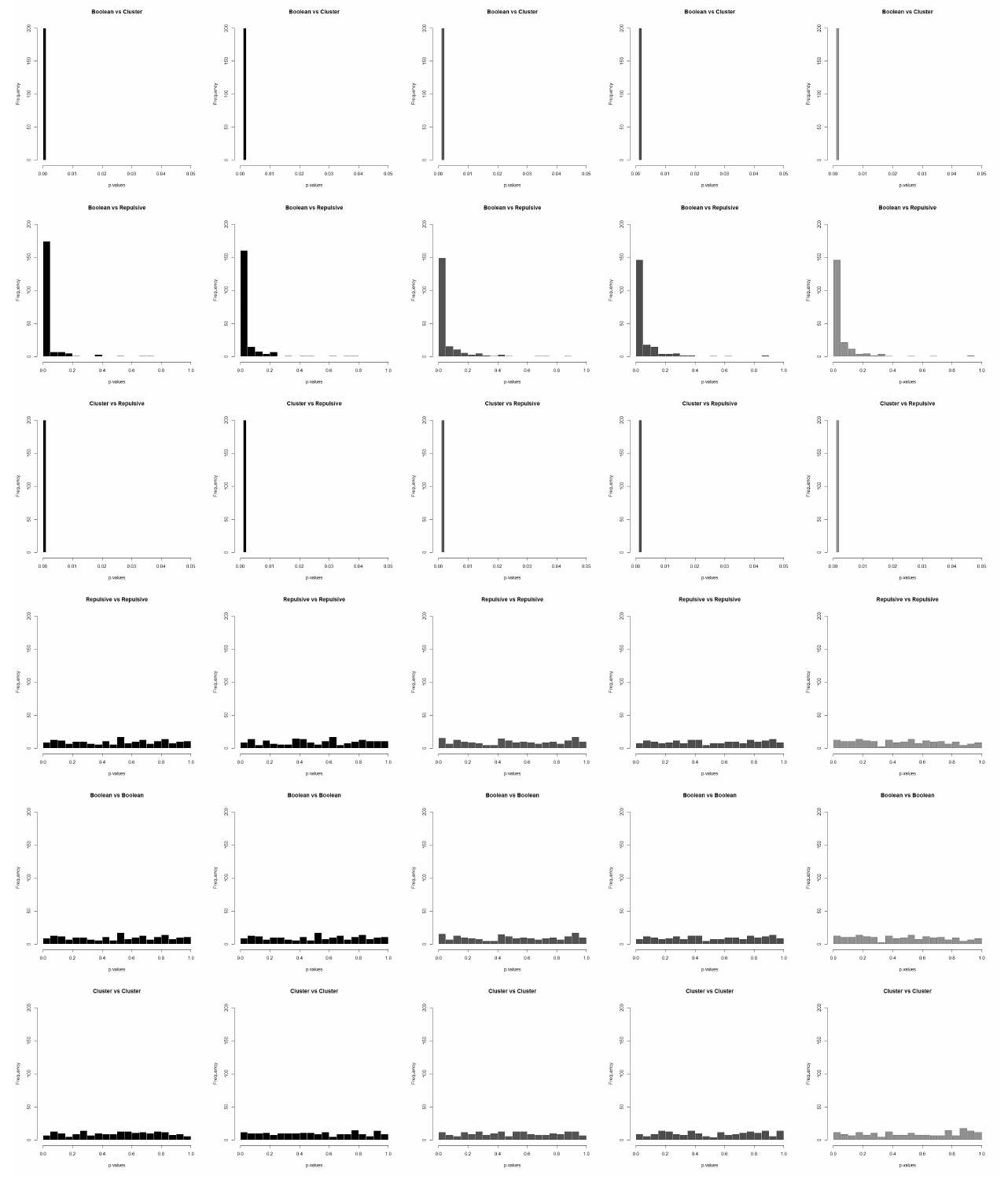}
\caption{Histograms of $p$-values of permutation version of test based on $\mathfrak{N}$-distances using kernel introduced in Example  \ref{ex:new_kernel1} (the first column), Example  \ref{ex:new_kernel2} (the second column), Cauchy kernel (the third column), Euclidean distance kernel (the fourth column) and Gaussian kernel (the fifth column) from Example \ref{ex:R_d_kernels} for testing all pairs of simulated processes,
   i.e. Boolean-cluster (the first row), Boolean-repulsive (the second row),
    cluster-repulsive (the third row), Boolean-Boolean (the forth row),
    cluster-cluster (the fifth row) and repulsive-repulsive (the sixth row).
    In each case, 200 pairs of realisations are used with number of sampled cells $m=100$, the number of permutations is $s=999.$}
\label{fig:histBCR_moje}
\end{figure}

\begin{figure}[t]
    \centering
\includegraphics[width=15cm]{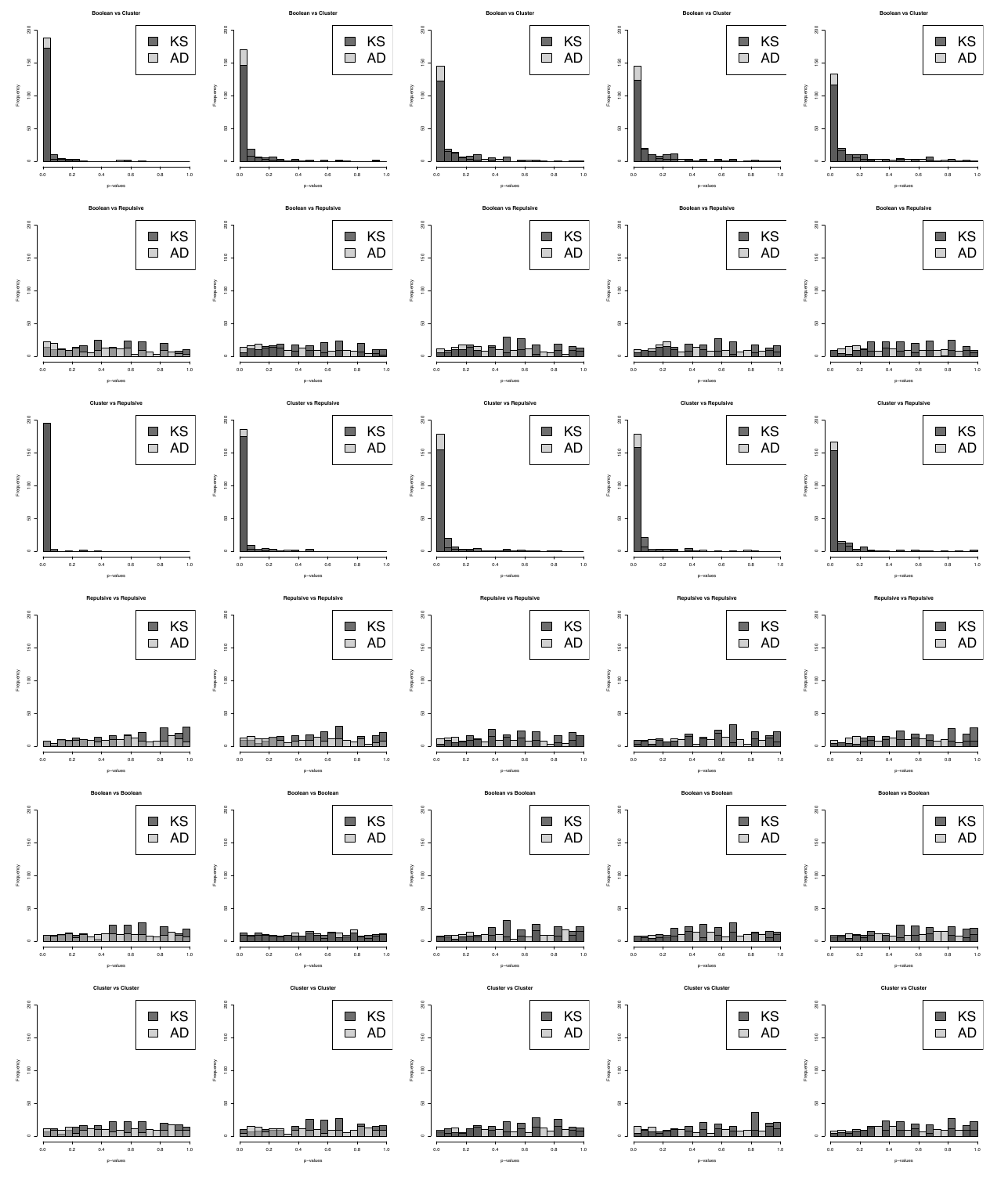}
\caption{Histograms of $p$-values when using version of test based on $\mathfrak{N}$-distances where KS and AD are involved using kernel introduced in Example \ref{ex:new_kernel1} (the first column), Example  \ref{ex:new_kernel2} (the second column), Cauchy kernel (the third column), Euclidean distance kernel (the fourth column) and Gaussian kernel (the fifth column) from Example \ref{ex:R_d_kernels} for testing all pairs of simulated processes,
   i.e. Boolean-cluster (upper left), Boolean-repulsive (upper middle),
    cluster-repulsive (upper right), Boolean-Boolean (lower left),
    cluster-cluster (lower middle) and repulsive-repulsive (lower right).
    In each case, 200 pairs of realisations are used with number of sampled cells $m=300$.}
\label{fig:histBCR_ks_ad}
\end{figure}

\begin{figure}[t]
\centering
\includegraphics[width=15cm]{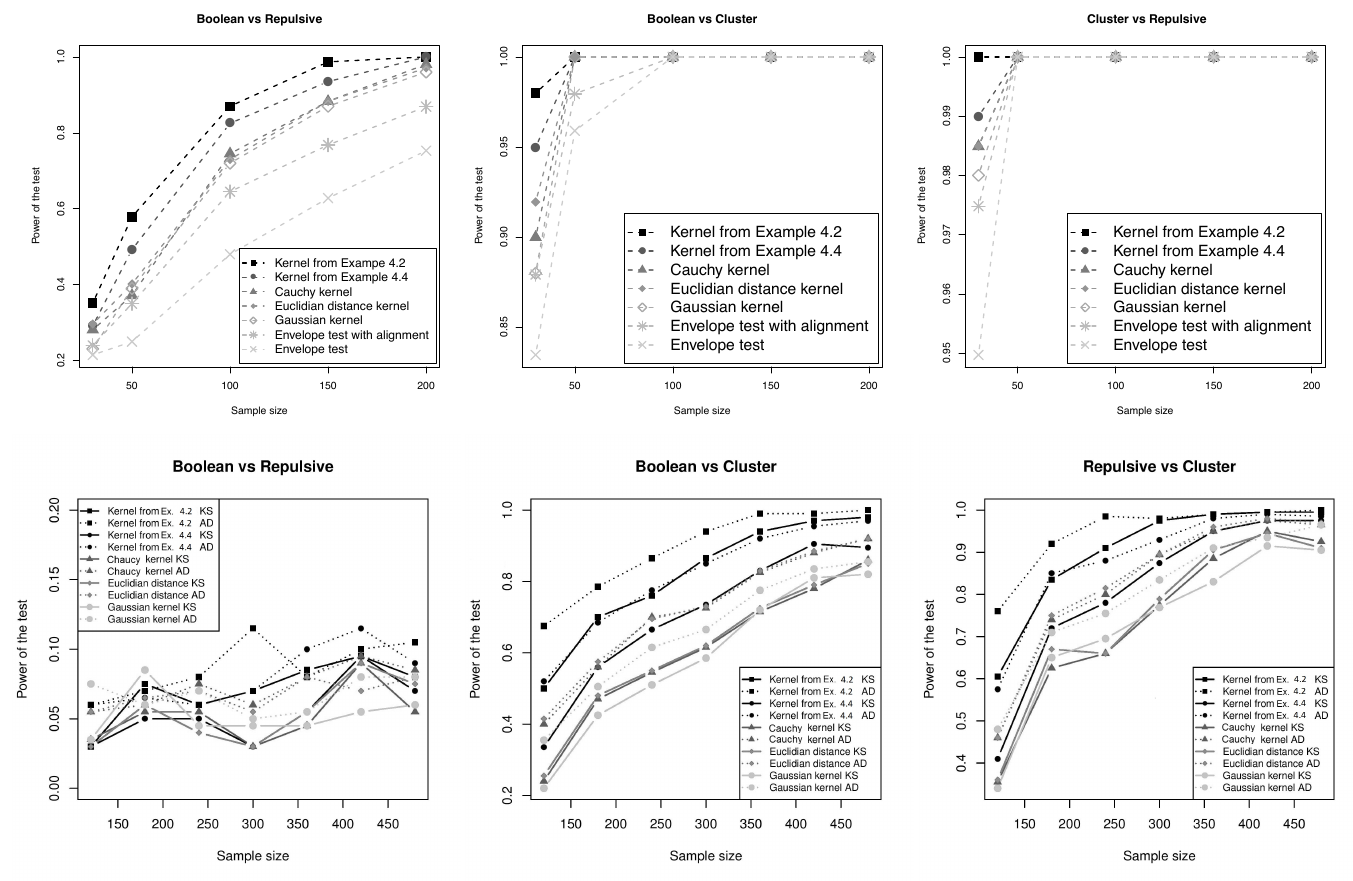}
\caption{Power of the tests for different cell sample sizes $m$ and employing different tests when comparing Boolean vs repulsive model (left), Boolean vs cluster model (middle) and cluster vs repulsive model (left). The significance level was chosen $\alpha=0.05.$ We compare the permutation versions of the tests based on $\mathfrak{N}$-distances and envelope tests (upper) and tests based on $\mathfrak{N}$-distances and KS and AD test (lower).}
\label{fig:power_env_Ntest}
\end{figure}

From the approximating tessellations, we randomly sample $m$ non-neighbouring cells 
whose support functions are evaluated in the angles
$u_i=\frac{2\pi}{n}i$, $i=0, \ldots, n-1$ for a given $n \in \mathbb{N}$, 
and we apply tests 
from Section \ref{sec:test} using different negative definite kernels on $\mathbb{R}^n \times \mathbb{R}^n$ applied on $n$-dimensional vectors of discretised support functions. 
Here, $n=10$ because computation of kernels are time consuming and wider simulation study showed that larger $n$ does not give any improvement in the test powers. 
The size of the sample of cells is $m=100$ for the permutation test while we use $s=999$ permutations, and $m=300$ for the test involving KS and AD.
As strongly negative definite kernels, we use the Euclidean distance kernel with $r=1,$ the Gaussian kernel with $ V=10 I$ ($ I$ is the identity matrix) and the Cauchy kernel with $V=I$ from Example $\ref{ex:R_d_kernels}$,
further the kernel from Example \ref{ex:new_kernel1}, where $D=3$ and $w_i=1,$
and the kernel from Example \ref{ex:new_kernel2}, where $D=3$ and $r=1.$
The histograms of 
$p$-values obtained from permutation 
test 
are plotted in Figure \ref{fig:histBCR_moje},
and $p$-values for the test using KS and AD are plotted in Figure \ref{fig:histBCR_ks_ad}.
In all cases, we can observe that comparing 
realisations of the same processes, 
the $p$-values are approximately uniformly distributed. 
It means that
the test finds the cells to be identically distributed, 
so the realisations are similar enough. 
On the other hand, when assessing similarity of different processes,
the $p$-values are very close to zero, so the hypothesis of identical
distribution of the cells is significantly often rejected, thus we can 
conclude that the similarity measure is very low, 
i.e. the original sets are not similar.
Finally, in Figure \ref{fig:power_env_Ntest},
it is seen that the new kernels give larger test power in both test procedures. 
Moreover, all permutation versions of kernel tests have larger power than envelope tests used in \cite{Getal16}, while the versions of kernel tests using KS and AD test show less power than envelope tests.

\section{Conclusion}
\label{sec:concl}

We have constructed the class of $\mathfrak{N}$-distances on the space of the distribution of random convex compact sets represented by random continuous functions on the unit sphere.
The obtained $\mathfrak{N}$-distances 
are used for testing the equality in distribution of random convex compact sets. 
Two different testing procedures have been proposed. 
Further,  the method for evaluating a truncated version of the strongly negative definite kernel used in the testing procedure is introduced, which allows to choose the depth of dependence structure.
A large rank of testing density functions is proposed.
Moreover, when observed the sets from different directions, 
different weights can be given to the directions. 
All these options make the test very flexible, so it allows to test hypotheses sensitive to different aspects.
Finally, the extensive simulation study including comparison of test powers
was performed with very satisfactory results.
Thus, we can conclude that we have obtained a convenient test for equality of distributions of random convex compact sets.
Note that although it was primarily constructed for distinguish between random sets, 
there is a wider field of applications, for example
comparing random continuous functions on compact domains or time series, the $p$-values obtained from the two-realisation test can be used as a similarity measure between two realisations of random sets in several machine learning algorithms etc.
Possible applications and corresponding modifications of the introduced methods are topics of the future research.

\bigskip

\section*{Acknowledgements}
The research was supported by The Czech Science Foundation, project No.~GA16-03708S.

\end{document}